\documentclass[11pt]{article}
\usepackage{amsmath,amssymb,amsthm}
\usepackage[margin=1.1in]{geometry} 
\usepackage[colorlinks=true,urlcolor=blue]{hyperref}
\usepackage{tikz}

\newtheorem{theorem}{Theorem}
\newtheorem{lemma}[theorem]{Lemma}
\newtheorem{proposition}[theorem]{Proposition}

\theoremstyle{definition}
\newtheorem{definition}[theorem]{Definition}
\newtheorem{example}[theorem]{Example}
\newtheorem{exercise}{Exercise}
\theoremstyle{remark}
\newtheorem{remark}[theorem]{Remark}

\newcommand{\R}{\mathbb{R}}

\title{Slicing a Doughnut at Random}
\author{Steven Finch\\
MIT Sloan School of Management}
\date{\today}

\begin{document}
\maketitle

\begin{abstract}
Choose a random plane that hits a solid torus. How likely is the slice to be
a single disk, and how fat should the doughnut be to make this as likely as
possible? Answering requires saying precisely what a ``random plane'' is. We
explain this from scratch, first for random lines in the plane and then for
random planes in space, and then solve the problem with nothing more than
calculus. Along the way we meet a two-dimensional cousin of the problem (the
doughnut's shadow) and two tempting wrong turns. One wrong turn appeared in
print; it turns out to be a correct solution of the two-dimensional cousin.
The other was the author's own first attempt, which luckily stayed in a drawer.
\end{abstract}

%======================================================================
\section{A doughnut and a knife}
%======================================================================

Let $C$ be the circle of radius $r\ge 1$ in the $xy$-plane, centered at the
origin, and let $T$ be the set of points of $\R^3$ at distance at most $1$
from $C$. Thus $T$ is a solid torus (a doughnut) with \emph{spine radius} $r$
and \emph{tube radius}~$1$; by Pappus's theorem its volume is
$\pi\cdot 2\pi r$.

Slice $T$ with a plane $P$. Apart from tangencies, the slice $P\cap T$ has
one of three shapes:
\begin{itemize}
\item a \emph{disk} (a cut near the outer rim),
\item an \emph{annulus} (a cut flat through the middle, like a bagel sliced
for toasting), or
\item \emph{two disks} (a vertical cut through the hole).
\end{itemize}

In 2012 I proposed the following problem to the \emph{Monthly}:
\begin{quote}
Let $P$ be a uniform random plane satisfying $P\cap T\neq\varnothing$. What
is the probability that the cross-section $P\cap T$ is simply connected? For
what value of $r$ is this probability maximized?
\end{quote}
``Simply connected'' means connected with no holes, so only the disk
counts.\footnote{Some solvers counted the two-disk slice as simply connected,
since each piece is; see the editorial comment in~\cite{Boukharfane}. We
follow the usual convention that a simply connected space must, in
particular, be connected.}

Before publication the editors asked me to say exactly what ``uniform random
plane'' means. The printed version~\cite{Finch} spells out a recipe: pick a
distance from the origin uniformly between $0$ and $1+r$, pick a normal
vector independently and uniformly on the unit sphere, and condition on $P$
meeting~$T$. The editors were right to ask. Problems about ``random''
geometric objects are notoriously slippery. Bertrand~\cite{Bertrand} famously
gave three different, equally plausible answers to one question about random
chords of a circle. This problem has two very tempting wrong answers, and
both of them actually appear later (Section~\ref{sec:wrong}).

Our goal is to explain carefully, assuming only calculus, what a random line
and a random plane are, and then to prove the following.

\begin{theorem}[Main result]\label{thm:main}
If $P$ is a uniform random plane meeting $T$, then
\[
\Pr\{P\cap T\text{ is a disk}\}
\;=\;\frac{\tfrac{\sqrt{r^2-1}}{r}+r\arcsin\tfrac1r}{1+\tfrac{\pi r}{4}}.
\]
This probability is maximal at a unique $r^*\approx 1.11890$, characterized by
\[
\psi=\cos\psi\,\Bigl(\sin\psi+\frac{\pi}{2}\Bigr),\qquad \psi=\arcsin\frac1{r^*}
\]
and the maximum value is $2\cos\psi=2\sqrt{1-1/r^{*2}}\approx 0.89721$.
\end{theorem}

The plan is as follows:
\begin{itemize}
\item Section~\ref{sec:lines}: random lines in the plane, where every idea
can be seen in its simplest form.
\item Section~\ref{sec:shadow}: a two-dimensional version of the problem,
the \emph{doughnut's shadow}.
\item Section~\ref{sec:planes}: random planes in space.
\item Section~\ref{sec:slice}: the proof of Theorem~\ref{thm:main}.
\item Section~\ref{sec:wrong}: the two wrong turns.
\item Section~\ref{sec:bonus}: a bonus fact and some exercises.
\end{itemize}

%======================================================================
\section{Random lines in the plane}\label{sec:lines}
%======================================================================

\subsection{Naming a line}
Every line in the plane has a unit normal vector
$u_\theta=(\cos\theta,\sin\theta)$ and a signed distance $p$ from the origin:
\[
L(\theta,p)=\{X\in\R^2:\ X\cdot u_\theta=p\},\qquad 0\le\theta<2\pi,\quad p\in\R .
\]
Each line has exactly two names, since $L(\theta,p)=L(\theta+\pi,-p)$, with
angles taken mod $2\pi$. This double counting is harmless, because we will
only ever compute \emph{ratios}. If we insist on $p\ge 0$, every line gets a
single name (lines through the origin aside, which never matter).

\subsection{What should ``uniform'' mean?}
There is no uniform \emph{probability} on the set of all lines, for the same
reason there is no uniform probability on all of $\R$: there is too much room.
What we can hope for is a fair way to \emph{measure} sets of lines. A uniform
random line inside any bounded region is then obtained by normalizing.

What should ``fair'' mean? Here is the guiding principle.

\begin{quote}
\emph{Fairness principle.} The measure of a set of lines should not depend on
where we put the origin or which way we point the axes.
\end{quote}

In other words, if we slide or turn a set of lines rigidly, its measure should
not change. Here is a measure that passes the test.

\begin{definition}
For a set $\mathcal L$ of lines and each angle $\theta$, let
$\mathcal L_\theta=\{p\in\R:\ L(\theta,p)\in\mathcal L\}$, and let $\lambda$
denote length. Define
\[
\mu(\mathcal L)=\int_0^{2\pi}\lambda(\mathcal L_\theta)\,d\theta .
\]
In words, $\mu(\mathcal L)$ is the area of the set of names $(\theta,p)$ of
the lines in $\mathcal L$.
\end{definition}

\begin{proposition}\label{prop:inv2}
The measure $\mu$ is unchanged by rotations, translations and reflections.
\end{proposition}

\begin{proof}
Rotation about the origin by an angle $\alpha$ carries $L(\theta,p)$ to
$L(\theta+\alpha,p)$. It merely shifts the integration variable $\theta$
(mod $2\pi$), so the integral is unchanged.

Reflection in the $x$-axis carries $L(\theta,p)$ to $L(-\theta,p)$, and the
same reasoning applies.

Translation by a vector $v$ carries $L(\theta,p)$ to
$L(\theta,\,p+v\cdot u_\theta)$, because $X\cdot u_\theta=p$ if and only if
$(X+v)\cdot u_\theta=p+v\cdot u_\theta$. For each fixed $\theta$, the set
$\mathcal L_\theta$ is slid along the $p$-axis by the fixed amount
$v\cdot u_\theta$, which does not change its length. Integrating over
$\theta$, the measure is unchanged.

Every rigid motion is a combination of these three kinds of maps.
\end{proof}

\begin{remark}
Up to a constant factor, $\mu$ is the \emph{only} measure on lines with this
invariance; see~\cite{Santalo} or~\cite{KlainRota}. So ``uniform random line''
has only one sensible meaning, and Jaynes~\cite{Jaynes} used exactly this idea
to settle Bertrand's paradox.

Other recipes that sound natural are not fair. Choosing a random slope and a
random intercept treats vertical lines very differently from horizontal ones.
Drawing the line through two random points of a disk favors lines through the
middle of the disk.
\end{remark}

\subsection{A random line that hits a given set}

\begin{definition}
Let $K$ be a bounded set in the plane. A \emph{uniform random line hitting $K$}
is a random line $L$ such that
\[
\Pr\{L\in\mathcal E\}=\frac{\mu(\mathcal E)}{\mu(\text{lines hitting }K)}
\]
for every set $\mathcal E$ of lines hitting $K$.
\end{definition}

\begin{proposition}[The recipe]\label{prop:recipe2}
Suppose $K$ lies inside the disk of radius $R$ about the origin. Choose
$\theta$ uniformly in $[0,2\pi)$ and, independently, $p$ uniformly in
$[0,R]$. If $L(\theta,p)$ misses $K$, throw it away and start again. The
first line that hits $K$ is a uniform random line hitting $K$, whatever the
value of $R$.
\end{proposition}

\begin{proof}
Every line hitting $K$ passes within distance $R$ of the origin, so it has a
unique name $(\theta,p)$ in the rectangle $[0,2\pi)\times[0,R]$. The point
$(\theta,p)$ is uniform on that rectangle. Throwing away failures leaves a
point that is uniform on the region of names of lines hitting $K$.
Probabilities are therefore ratios of areas of sets of names, which are
ratios of $\mu$-measures. (Using $p\ge 0$ halves both areas, which does not
affect the ratio.)
\end{proof}

\begin{example}[Bertrand's chord]
A uniform random line hitting the unit circle cuts a chord longer than
$\sqrt3$ (the side of the inscribed equilateral triangle) exactly when
$p<\tfrac12$. With $R=1$ in the recipe, $p$ is uniform on $[0,1]$, so the
probability is $\tfrac12$.
\end{example}

\subsection{Warning: conditioning changes the odds of directions}\label{sec:stick}

Let $K$ be a stick, the segment from $(-1,0)$ to $(1,0)$. The signed
distances of its endpoints from $L(\theta,p)$ are $\pm\cos\theta-p$, so the
line hits the stick exactly when $|p|\le|\cos\theta|$.

In the recipe with $R=1$, a line with normal angle $\theta$ therefore
survives with probability $|\cos\theta|$. So, \emph{given that the line hits
the stick}, its normal angle is \emph{not} uniform: its density is
$|\cos\theta|/4$ on $[0,2\pi)$. Lines that hit a stick prefer to cross it
broadside. This is no paradox: a broadside line has more room to hit.

The moral is worth stating in symbols. For any event $\mathcal E$,
\[
\Pr\{\mathcal E\mid L\text{ hits }K\}
=\frac{\displaystyle\int\lambda(\mathcal E_\theta)\,d\theta}
{\displaystyle\int\lambda(\text{hits}_\theta)\,d\theta}
\qquad\text{(a ratio of two totals)},
\]
\emph{not} the average over $\theta$ of the ratios
$\lambda(\mathcal E_\theta)/\lambda(\text{hits}_\theta)$. The second
expression secretly pretends that the direction is still uniform after
conditioning.

%======================================================================
\section{Warm-up: the doughnut's shadow}\label{sec:shadow}
%======================================================================

We use coordinates $(x,z)$ in the plane, because in Section~\ref{sec:slice}
this plane will be the $xz$-plane in space. Let $D_+$ and $D_-$ be unit disks
(``coins'') centered at $(r,0)$ and $(-r,0)$, where $r\ge 1$. Let $R$ be the
\emph{racetrack}: the set of points within distance $1$ of the segment joining
the two centers. Equivalently, $R$ is the convex hull of the two coins:
\begin{equation}\label{eq:R}
R=\bigl\{(x,z):\ |z|\le 1,\ |x|\le r+\sqrt{1-z^2}\bigr\}.
\end{equation}

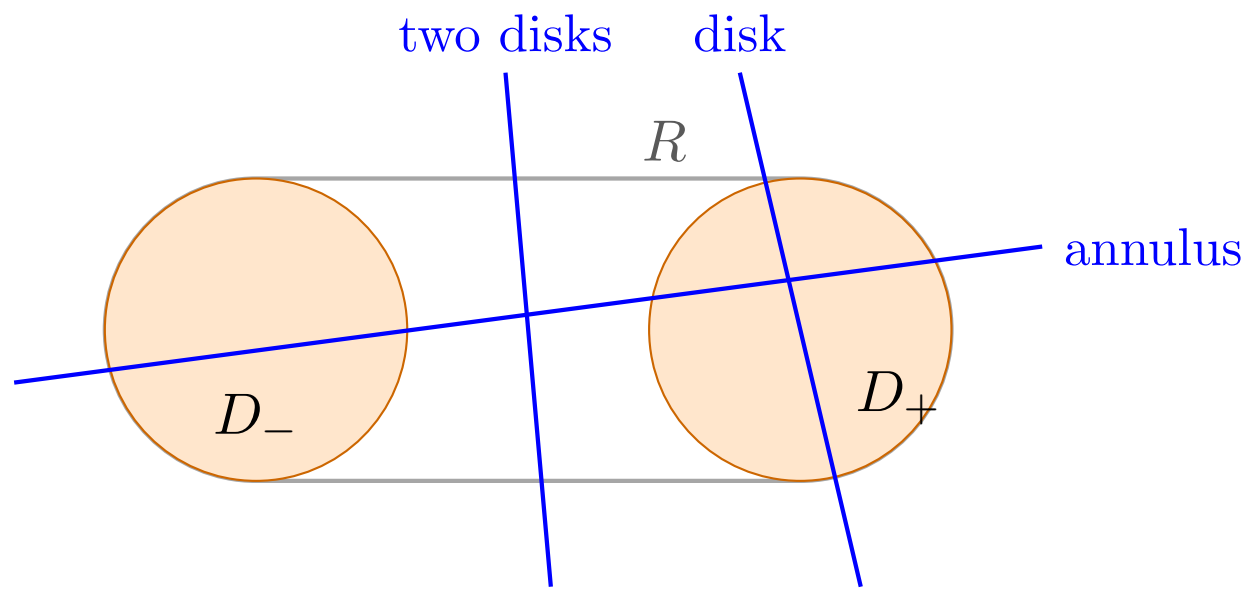
\begin{figure}[ht]
\centering
\begin{tikzpicture}[scale=1.0]
\def\rr{1.8}
\draw[thick,gray!70] (-\rr,1) -- (\rr,1) arc (90:-90:1) -- (-\rr,-1) arc (270:90:1) -- cycle;
\filldraw[fill=orange!20,draw=orange!80!black] (\rr,0) circle (1);
\filldraw[fill=orange!20,draw=orange!80!black] (-\rr,0) circle (1);
\node at (2.45,-0.45) {$D_+$};
\node at (-\rr,-0.6) {$D_-$};
\node[gray!70!black] at (0.9,1.25) {$R$};
\draw[blue,thick] (2.2,-1.7) -- (1.4,1.7) node[above] {\small disk};
\draw[blue,thick] (-3.4,-0.35) -- (3.4,0.55) node[right] {\small annulus};
\draw[blue,thick] (0.15,-1.7) -- (-0.15,1.7) node[above] {\small two disks};
\end{tikzpicture}
\caption{The doughnut's shadow, drawn with $r=1.8$: two coins $D_\pm$ inside
the racetrack $R$. A line meeting $R$ may hit one coin, both coins or neither.
In Section~\ref{sec:slice} these three cases become the three kinds of slices
of the torus.}
\label{fig:shadow}
\end{figure}

Why this picture? Lay the doughnut $T$ on a table and look at it from the
side, along the $y$-axis. Its outline is the racetrack $R$, and if you cut it
along the $xz$-plane you see the two coins. We will see in
Section~\ref{sec:slice} that the torus problem is \emph{literally} the
following problem, viewed from the side.

\begin{quote}
\emph{Shadow problem.} A uniform random line hits the racetrack $R$. What is
the probability that it hits exactly one coin?
\end{quote}

Write each line as
$\ell(\beta,p)=\{(x,z):\ x\sin\beta+z\cos\beta=p\}$, with $0\le\beta<2\pi$ and
$p\ge0$. This is the naming of Section~\ref{sec:lines} with
$\theta=\tfrac{\pi}{2}-\beta$, so $\mu$ is still ``$d\beta\,dp$''. Put
\[
c=r\,|\sin\beta| .
\]
Since the coin centers are $(\pm r,0)$, the number $c$ is how far the centers
stick out in the direction of the normal.

\begin{lemma}\label{lem:shadow}
For $p\ge 0$:
\begin{enumerate}
\item[(i)] $\ell(\beta,p)$ meets $R$ if and only if $p\le c+1$;
\item[(ii)] $\ell(\beta,p)$ meets exactly one coin if and only if
$|c-p|\le 1<c+p$.
\end{enumerate}
Consequently, for each fixed $\beta$, the set of $p\ge0$ for which the line
meets $R$ has length $c+1$, and the set for which it meets exactly one coin
has length $2\min(c,1)$.
\end{lemma}

\begin{proof}
By the symmetry $x\mapsto -x$ we may assume $\sin\beta\ge0$. The signed
distance from a point $(x_0,z_0)$ to the line is
$x_0\sin\beta+z_0\cos\beta-p$. For the two centers this gives $c-p$ and
$-c-p$.

(i) Along the segment joining the centers, the signed distance varies
linearly from $-c-p$ to $c-p$. If $p\le c$, the line crosses the segment and
so meets $R$. If $p>c$, every point of the segment is at distance at least
$p-c$ from the line, with equality at $(r,0)$. So the line meets $R$
exactly when $p-c\le 1$.

(ii) The line meets $D_+$ if and only if $|c-p|\le1$, and meets $D_-$ if and
only if $c+p\le 1$. Since $|c-p|\le c+p$, hitting $D_-$ forces hitting $D_+$.
So ``exactly one'' means $|c-p|\le 1<c+p$.

For the lengths:
\begin{itemize}
\item if $c\ge 1$, the condition in (ii) is $p\in[c-1,c+1]$, of length $2$;
\item if $c<1$, it is $p\in(1-c,1+c]$, of length $2c$.
\end{itemize}
In both cases the length is $2\min(c,1)$.
\end{proof}

\begin{theorem}[The shadow problem]\label{thm:shadow}
The probability that a uniform random line hitting $R$ meets exactly one coin
is
\[
P_2(r)=\frac{\displaystyle\int_0^{\pi/2}2\min(r\sin\beta,1)\,d\beta}
{\displaystyle\int_0^{\pi/2}(1+r\sin\beta)\,d\beta}
=\frac{2r+2\arccos(1/r)-2\sqrt{r^2-1}}{r+\pi/2}.
\]
It is maximal at a unique $r\approx1.24376$. There
$\varphi=\arccos(1/r)$ solves $\varphi=\frac{\pi}{2}(1-\sin\varphi)$, and the
maximum value is $4\varphi/\pi\approx 0.81078$.
\end{theorem}

\begin{proof}
By Lemma~\ref{lem:shadow} and the moral of Section~\ref{sec:stick}, the
probability is the ratio of the integrals of $2\min(c,1)$ and $c+1$ over
$\beta\in[0,2\pi)$. Since $c$ depends only on $|\sin\beta|$, each integral is
four times the integral over $[0,\pi/2]$.

The denominator is $\pi/2+r$. For the numerator, let
$\beta_0=\arcsin(1/r)$. Then
\[
\int_0^{\beta_0}2r\sin\beta\,d\beta+\int_{\beta_0}^{\pi/2}2\,d\beta
=2r-2r\cos\beta_0+\pi-2\beta_0
=2r-2\sqrt{r^2-1}+2\arccos(1/r).
\]

To find the maximum, write $N(r)$ for the numerator and $D(r)=r+\pi/2$, and
let $\varphi=\arccos(1/r)$, so that $\sqrt{r^2-1}=r\sin\varphi$. Then
$N=2r(1-\sin\varphi)+2\varphi$, and a short computation gives
$N'(r)=2-2\sqrt{r^2-1}/r=2(1-\sin\varphi)$. Hence
\[
N'D-ND'=\pi(1-\sin\varphi)-2\varphi=-2G(\varphi),\qquad
G(\varphi)=\varphi-\tfrac{\pi}{2}(1-\sin\varphi).
\]
The function $G$ is increasing, with $G(0)=-\frac{\pi}{2}$ and
$G(\frac{\pi}{2})=\frac{\pi}{2}$, and $\varphi$ increases with $r$. So $P_2$
increases and then decreases, with a single maximum where $G(\varphi)=0$.
There $P_2=N'/D'=2(1-\sin\varphi)=4\varphi/\pi$. Newton's method gives
$\varphi\approx 0.63678$ and $r=1/\cos\varphi\approx 1.24376$.
\end{proof}

\begin{remark}[A Crofton check]
Counting each line once, the measure of lines hitting $R$ is
$\int_0^{2\pi}(1+r|\sin\beta|)\,d\beta=2\pi+4r$, which is exactly the
perimeter of the racetrack. This is no accident. \emph{Crofton's formula}
says that the measure of the set of lines hitting a convex set equals its
perimeter~\cite{Santalo,KlainRota}.
\end{remark}

\begin{remark}[Other flat doughnuts]
The most obvious flat doughnut is the annulus of points within distance $1$
of a circle of radius $r$. A random line hitting it meets it in one piece
with probability $2/(1+r)$ (Exercise~\ref{ex:annulus}). This is largest at
$r=1$, so there is nothing to find. The annulus is too round to be
interesting: all directions look alike. The racetrack is the honest analog,
because it really is the side view of the torus.
\end{remark}

%======================================================================
\section{Random planes in space}\label{sec:planes}
%======================================================================

Everything carries over, with one new twist: choosing a random
\emph{direction} in space.

A plane is named by a unit normal $n\in S^2$ (the unit sphere) and a signed
distance $m$:
\[
\Pi(n,m)=\{X\in\R^3:\ X\cdot n=m\}.
\]
As before, $\Pi(n,m)=\Pi(-n,-m)$, so every plane has two names, or one name
if we insist on $m\ge 0$.

\begin{definition}
For a set $\mathcal P$ of planes, let
$\mathcal P_n=\{m:\ \Pi(n,m)\in\mathcal P\}$ and
\[
\mu(\mathcal P)=\int_{S^2}\lambda(\mathcal P_n)\,dA(n),
\]
where $dA$ is surface area on the unit sphere. In words: measure directions
by \emph{area on the sphere} and distances by \emph{length}.
\end{definition}

\begin{proposition}
The measure $\mu$ is unchanged by rotations, translations and reflections.
\end{proposition}

\begin{proof}
A rotation or reflection $Q$ fixing the origin carries $\Pi(n,m)$ to
$\Pi(Qn,m)$, and $Q$ preserves area on the sphere. A translation by $v$
carries $\Pi(n,m)$ to $\Pi(n,m+v\cdot n)$, which slides each set
$\mathcal P_n$ without changing its length.
\end{proof}

As in the plane, $\mu$ is the only invariant measure on planes, up to a
constant factor~\cite{Santalo}. This is what ``uniform random plane'' means.

\subsection{Choosing a direction: the latitude trap}
How do we pick $n$ uniformly on the sphere, that is, with probabilities
proportional to area? It is tempting to pick the angle $\beta$ between $n$
and the $z$-axis uniformly in $[0,\pi]$. That is wrong: it crowds directions
near the poles, where circles of latitude are short. The right answer goes
back to Archimedes.

\begin{lemma}[Archimedes' hat-box theorem]\label{lem:archimedes}
For $-1\le a\le b\le 1$, the area of the zone\\ $\{n\in S^2:\ a\le n_z\le b\}$
is $2\pi(b-a)$.
\end{lemma}

\begin{proof}
The sphere is the surface of revolution of $f(z)=\sqrt{1-z^2}$ about the
$z$-axis. Its area between heights $a$ and $b$ is
$\int_a^b 2\pi f(z)\sqrt{1+f'(z)^2}\,dz$. Here
$1+f'(z)^2=1/(1-z^2)$, so the integrand is the constant $2\pi$.
\end{proof}

So for a uniform direction, the height $n_z=\cos\beta$ is uniform on
$[-1,1]$. Equivalently, $\beta$ has density $\frac12\sin\beta$, not a
constant. Cutting the sphere into thin zones and applying
Lemma~\ref{lem:archimedes} to each one gives, for any reasonable function
$g$,
\begin{equation}\label{eq:arch}
\int_{S^2}g(|n_z|)\,dA(n)=4\pi\int_0^1 g(u)\,du .
\end{equation}

\subsection{A random plane that hits a given body}
Exactly as in Proposition~\ref{prop:recipe2}, define a \emph{uniform random
plane hitting $T$} by $\Pr\{\mathcal E\}=\mu(\mathcal E)/\mu(\text{planes
hitting }T)$. Here is a recipe for producing one.
\begin{itemize}
\item Choose $n$ uniformly on $S^2$ (for instance with $n_z$ uniform on
$[-1,1]$ and a uniform longitude).
\item Choose $m$ uniformly on $[0,1+r]$, independently of $n$.
\item Throw the plane away unless it meets $T$, and start again.
\end{itemize}
This is precisely the recipe in the published wording of the
problem~\cite{Finch}. Because the fair measure is unique, the original phrase
``uniform random plane'' meant the same thing. And just as with the stick,
\emph{given that $P$ meets $T$, its normal is no longer uniform}.

%======================================================================
\section{Slicing the doughnut}\label{sec:slice}
%======================================================================

The nearest point of $C$ to $(x,y,z)$ lies in the direction of $(x,y,0)$, so
the distance from $(x,y,z)$ to $C$ is
$\sqrt{(\sqrt{x^2+y^2}-r)^2+z^2}$. Hence
\[
T=\Bigl\{(x,y,z):\ \bigl(\sqrt{x^2+y^2}-r\bigr)^2+z^2\le 1\Bigr\}.
\]

\paragraph{Normalizing the plane.}
The torus is unchanged by rotations about the $z$-axis and by the reflection
$z\mapsto -z$. Using these, together with the $m\ge0$ naming, we may move any
plane to one of the form $\Pi(n,m)$ with
\[
n=(\sin\beta,0,\cos\beta),\qquad 0\le\beta\le\tfrac{\pi}{2},\qquad m\ge 0,
\]
without changing the shape of the slice. For the original normal,
$\sin\beta=\sqrt{n_x^2+n_y^2}$ and $\cos\beta=|n_z|$. Such a plane contains
the $y$-direction, so it stands perpendicular to the $xz$-plane and meets it
in the line
\[
\ell:\ x\sin\beta+z\cos\beta=m .
\]
This is the line $\ell(\beta,m)$ of Section~\ref{sec:shadow}. The $xz$-plane
cuts $T$ in the two coins $D_\pm$, and seen from the side, $P$ looks like
$\ell$. Put $c=r\sin\beta$, as before.

\begin{lemma}\label{lem:meets}
$P$ meets $T$ if and only if $\ell$ meets the racetrack $R$, that is, if and
only if $m\le c+1$.
\end{lemma}

\begin{proof}
The equation of $P$ does not involve $y$. So $P$ meets $T$ exactly when some
point of $T$ has its shadow $(x,z)$ on $\ell$, that is, when $\ell$ meets the
shadow of $T$ on the $xz$-plane.

That shadow is $R$. Indeed, $(x,z)$ is a shadow point when some $y$ satisfies
$|z|\le1$ and $r-w\le\sqrt{x^2+y^2}\le r+w$, where $w=\sqrt{1-z^2}$. As $y$
varies, $\sqrt{x^2+y^2}$ takes every value in $[|x|,\infty)$, so the
condition is $|x|\le r+w$. This is \eqref{eq:R}.

The last statement is Lemma~\ref{lem:shadow}(i).
\end{proof}

\begin{lemma}[Slicing lemma]\label{lem:slice}
Suppose $P$ meets $T$ and $\ell$ is not tangent to either coin. Then $P\cap
T$ is
\begin{enumerate}
\item[(a)] a disk if $\ell$ meets exactly one coin;
\item[(b)] an annulus if $\ell$ meets both coins;
\item[(c)] two disjoint disks if $\ell$ meets neither coin.
\end{enumerate}
\end{lemma}

\begin{proof}
Let $e=(\cos\beta,0,-\sin\beta)$ and $j=(0,1,0)$. These are perpendicular unit
vectors parallel to $P$, so every point of $P$ is
$X(s,y)=mn+se+yj$ for a unique pair $(s,y)$, and $(s,y)$ are Cartesian
coordinates on $P$. The points with $y=0$ form $\ell$; write
$q(s)=mn+se=(x(s),0,z(s))$.

\smallskip
\emph{Step 1: vertical fibers.} Fix $s$ and write $x=x(s)$, $z=z(s)$. The
point $(x,y,z)$ lies in $T$ when $|z|\le 1$ and
$r-w\le\sqrt{x^2+y^2}\le r+w$, where $w=\sqrt{1-z^2}\in[0,1]$. As $|y|$ grows
from $0$, $\sqrt{x^2+y^2}$ increases from $|x|$. So the allowed values of $y$
form the set $a(s)\le|y|\le b(s)$, where
\[
b(s)=\sqrt{(r+w)^2-x^2},\qquad a(s)=\sqrt{\max\{0,\,(r-w)^2-x^2\}}.
\]
Two observations follow.
\begin{itemize}
\item The fiber is nonempty exactly when $q(s)\in R$, by the proof of
Lemma~\ref{lem:meets}. Since $R$ is convex, these $s$ form a closed interval
$[s_1,s_2]$. In its interior $a(s)<b(s)$, because there $|z|<1$ and
$b^2-a^2=\min\{(r+w)^2-x^2,\,4rw\}>0$.
\item $a(s)=0$ exactly when $r-w\le|x|\le r+w$, that is, when
$(|x|-r)^2+z^2\le1$, that is, when $q(s)$ lies in a coin.
\end{itemize}

\smallskip
\emph{Step 2: two halves.} Let
$U=\{(s,y):\ s_1\le s\le s_2,\ a(s)\le y\le b(s)\}$ be the part of the slice
with $y\ge 0$, and let $L$ be its mirror image under $y\mapsto -y$. The set
$U$ lies between the graphs of two continuous functions that pinch together
only at the ends, so it is a closed disk (think of a squashed pancake), and
so is $L$. The slice is $U\cup L$, and
\[
U\cap L=\{(s,0):\ a(s)=0\}=\ell\cap(D_+\cup D_-),
\]
which consists of zero, one or two disjoint chords. These lie along the
bottom edge of $U$ and the top edge of $L$.

\smallskip
\emph{Step 3: gluing.} Two pancakes glued along nothing are two disks. Glued
along one seam, they form a single disk. Glued along two separate seams, they
enclose a hole between the seams and form an annulus.
\end{proof}

Combining Lemmas~\ref{lem:shadow}, \ref{lem:meets} and~\ref{lem:slice}: for
the normalized plane,
\begin{equation}\label{eq:criteria}
P\text{ meets }T\iff m\le c+1,\qquad
P\cap T\text{ is a disk}\iff |c-m|\le 1<c+m,
\end{equation}
with $c=r\sin\beta=r\sqrt{1-n_z^2}$. These are \emph{word for word} the
conditions of the shadow problem.

\begin{proof}[Proof of Theorem~\ref{thm:main}]
Name planes with $m\ge0$. By \eqref{eq:criteria} and
Lemma~\ref{lem:shadow}, for each normal $n$ the set of $m$ giving a disk has
length $2\min(c,1)$, and the set of $m$ giving a hit has length $1+c$. Both
depend only on $|n_z|$. By \eqref{eq:arch}, with $u=|n_z|$,
\[
\Pr\{\text{disk}\}=\frac{\int_{S^2}2\min(c,1)\,dA}{\int_{S^2}(1+c)\,dA}
=\frac{\int_0^1 2\min\bigl(r\sqrt{1-u^2},1\bigr)\,du}{\int_0^1\bigl(1+r\sqrt{1-u^2}\bigr)\,du}.
\]

The denominator is $1+\pi r/4$, since $\int_0^1\sqrt{1-u^2}\,du$ is the area
of a quarter of the unit disk.

For the numerator, $r\sqrt{1-u^2}\ge1$ exactly when $u\le
s:=\sqrt{1-1/r^2}$. Using
$\int\sqrt{1-u^2}\,du=\frac12\bigl(u\sqrt{1-u^2}+\arcsin u\bigr)$ and
$\sqrt{1-s^2}=1/r$,
\[
2s+2r\int_s^1\sqrt{1-u^2}\,du
=2s+2r\Bigl(\frac{\pi}{4}-\frac{s}{2r}-\frac{\arcsin s}{2}\Bigr)
=s+r\arccos s .
\]
Since $s=\sqrt{r^2-1}/r$ and $\arccos s=\arcsin(1/r)$, this proves the
formula.

For the maximum, write $N(r)=\sqrt{r^2-1}/r+r\arcsin(1/r)$ and
$D(r)=1+\pi r/4$. Differentiating,
$N'(r)=\arcsin(1/r)-\sqrt{r^2-1}/r^2$. Let $\psi=\arcsin(1/r)$, so that
$\sin\psi=1/r$, $\cos\psi=\sqrt{r^2-1}/r$ and
$\sqrt{r^2-1}/r^2=\sin\psi\cos\psi$. Then
\[
N=\cos\psi+\frac{\psi}{\sin\psi},\qquad N'=\psi-\sin\psi\cos\psi,\qquad
D=1+\frac{\pi}{4\sin\psi},\qquad D'=\frac{\pi}{4},
\]
and the terms $\pm\pi\psi/(4\sin\psi)$ cancel in
\[
N'D-ND'=F(\psi):=\psi-\sin\psi\cos\psi-\frac{\pi}{2}\cos\psi .
\]
So $P'(r)=F(\psi)/D^2$.

Now $F'(\psi)=2\sin^2\psi+\frac{\pi}{2}\sin\psi>0$, while
$F(0)=-\frac{\pi}{2}$ and $F(\frac{\pi}{2})=\frac{\pi}{2}$. Also $\psi$
decreases from $\frac{\pi}{2}$ to $0$ as $r$ runs from $1$ to $\infty$. Hence
$P$ increases and then decreases, with a unique maximum where $F(\psi)=0$,
that is, $\psi=\cos\psi\,(\sin\psi+\frac{\pi}{2})$.

At that point $N/D=N'/D'=\frac{4}{\pi}(\psi-\sin\psi\cos\psi)=2\cos\psi$.
Newton's method gives $\psi\approx 1.105595$, so
$r^*=1/\sin\psi\approx 1.11890$ and $P\approx 0.89721$.
\end{proof}

\paragraph{The moral, in one display.}
In terms of the tilt angle $\beta$,
\[
P_2(r)=\frac{\int_0^{\pi/2}2\min(r\sin\beta,1)\,d\beta}{\int_0^{\pi/2}(1+r\sin\beta)\,d\beta},
\qquad
P(r)=\frac{\int_0^{\pi/2}2\min(r\sin\beta,1)\,\sin\beta\,d\beta}{\int_0^{\pi/2}(1+r\sin\beta)\,\sin\beta\,d\beta}.
\]
\emph{The torus problem is the shadow problem with the directions reweighted
by Archimedes.} Lines in the plane have uniformly distributed directions;
normals in space do not.

%======================================================================
\section{Two tempting wrong turns}\label{sec:wrong}
%======================================================================

\paragraph{Wrong turn 1: a uniform tilt angle.}
If one chooses the tilt $\beta$ uniformly (the latitude trap), the
computation goes through as above with weight $1$ in place of $\sin\beta$.
The result is $P_2(r)$ of Theorem~\ref{thm:shadow}, with maximum $\approx
0.81078$ at $r\approx 1.24376$. This is the formula printed
in~\cite{Boukharfane}.\footnote{There the critical equation is printed as
$\sin\bigl(\frac{\pi}{2}\sqrt{1-(r/2)^2}\bigr)=\frac1r$; the expression
under the radical should be $1-(1/r)^2$. The corrected equation is
equivalent to $\varphi=\frac{\pi}{2}(1-\sin\varphi)$ with
$\cos\varphi=1/r$.} As we now know, it is a perfectly correct answer, 
but to the shadow problem!

\paragraph{Wrong turn 2: averaging ratios.}
For a fixed direction, the conditional probability of a disk given a hit is
$2\min(c,1)/(1+c)$. It is tempting to average this against the (correct!)
direction density $\sin\beta$:
\[
q(r)=\int_0^{\pi/2}\frac{2\min(r\sin\beta,1)}{1+r\sin\beta}\,\sin\beta\,d\beta .
\]
This was my own first answer in 2012. It gives, for instance,
$q(1)=4-\pi\approx 0.8584$, with a maximum $\approx 0.8781$ at $r\approx
1.14128$.

The flaw is the stick lesson of Section~\ref{sec:stick}. Given that $P$
meets $T$, directions in which the doughnut looks wide are more likely: the
conditional density of $\beta$ is
$(1+r\sin\beta)\sin\beta/(1+\pi r/4)$, not $\sin\beta$. Averaging the
per-direction ratio against this corrected density makes the factors
$1+r\sin\beta$ cancel, and gives back Theorem~\ref{thm:main}.

\begin{center}
\renewcommand{\arraystretch}{1.25}
\begin{tabular}{lcccc}
 & weight on $\beta$ & value at $r=1$ & best $r$ & best value\\ \hline
Wrong turn 1 (= shadow problem) & $1$ & $0.7780$ & $1.24376$ & $0.81078$\\
Wrong turn 2 (average of ratios) & $\sin\beta$, unconditioned & $0.8584$ & $1.14128$ & $0.8781$\\
Correct (Theorem~\ref{thm:main}) & $\sin\beta$, conditioned & $0.8798$ & $1.11890$ & $0.89721$
\end{tabular}
\end{center}

%======================================================================
\section{A bonus, and some exercises}\label{sec:bonus}
%======================================================================

\paragraph{Bonus.}
Score a disk as $1$, two disks as $2$ and an annulus as $0$; this score is
the Euler characteristic $\chi$ of the slice. For each direction, the
integral of $\chi$ over $m\ge0$ is
$2\min(c,1)+2\max(c-1,0)=2c$, a pleasantly simple answer. By
\eqref{eq:arch},
\[
\mathbb E\bigl[\chi(P\cap T)\mid P\text{ meets }T\bigr]
=\frac{\pi r/2}{1+\pi r/4}=\frac{2\pi r}{4+\pi r}.
\]
At $r=1$ two-disk slices are impossible ($c\le 1$), so $\Pr\{\text{disk}\}$
must equal $2\pi/(4+\pi)\approx 0.8798$, and it does.

\begin{exercise}\label{ex:annulus}
Let $A$ be the annulus of points within distance $1$ of a circle of radius
$r$ in the plane. Show that a uniform random line hitting $A$ meets it in a
single segment with probability $2/(1+r)$.
\end{exercise}

\begin{exercise}
In the shadow picture, condition instead on the line hitting at least one
coin. Show that the probability of hitting exactly one is $4g/(\pi+2g)$,
where $g=r+\arccos(1/r)-\sqrt{r^2-1}$, and that this increases with $r$
towards $1$. Why is the racetrack, rather than the pair of coins, the right
analog of the torus?
\end{exercise}

\begin{exercise}
Show that the slice of $T$ is an annulus with probability
\[
\frac{1-\frac{\pi r}{4}+\frac r2\arccos\frac1r-\frac{\sqrt{r^2-1}}{2r}}{1+\frac{\pi r}{4}},
\]
and is two disks with probability
\[
\frac{\frac r2\arccos\frac1r-\frac{\sqrt{r^2-1}}{2r}}{1+\frac{\pi r}{4}}.
\]
Check that, together with Theorem~\ref{thm:main}, the three probabilities
add up to $1$.
\end{exercise}

\begin{exercise}
Simulate! Program the recipe of Section~\ref{sec:planes}, classify each slice
using \eqref{eq:criteria}, and estimate the probability of a disk at
$r=1.1189$. Then program the two wrong turns and watch all three estimates
disagree.
\end{exercise}

\paragraph{Acknowledgment.}
I thank the editors of the \emph{Monthly} for insisting, back in 2012, that I
specify what ``random'' means. They saved me from my own Wrong Turn~2.
Conversations with Anthropic Claude Opus 5.5 and OpenAI GPT-6 Luna were 
essential in the preparation of this paper.  No AI funding whatsoever for 
this study was received.

\end{document}